\documentclass[reqno]{amsart}
\usepackage{hyperref}
\usepackage{setspace}
\usepackage{graphicx}
\usepackage{amssymb}
\usepackage{mathrsfs}

\begin{document}
    \title[\hfilneg ]
    { EXISTENCE AND UNIQUENESS OF
        BOUNDARY VALUE PROBLEMS FOR
        HILFER-HADAMARD-TYPE FRACTIONAL
        DIFFERENTIAL EQUATIONS}

 \author[\hfil\hfilneg]{Ahmad Y. A. Salamooni, D. D. Pawar }
 \address{Ahmad Y. A. Salamooni \newline
     School of Mathematical Sciences, Swami Ramanand Teerth Marathwada University, Nanded-431606, India}
      \email{ayousss83@gmail.com}

 \address{D. D. Pawar \newline
    School of Mathematical Sciences, Swami Ramanand Teerth Marathwada University, Nanded-431606, India}
     \email{dypawar@yahoo.com}

     \keywords{Existence,uniqueness,boundary value problems,
    Hilfer-Hadamard type, fractional differential equation and fractional calculus}

    \begin{abstract}
    In this paper, we used some theorems of fixed point for studying the results of
    existence and uniqueness  for Hilfer-Hadamard-Type fractional differential equations,
     \[_{H}D^{\alpha,\beta}x(t)+f(t,x(t))=0, ~~~~~~        on~~the~~ interval~~  J:=(1,e]\]
    with boundary value problems
   \[x(1+\epsilon)=0,~~~~~~~~~~_{H}D^{1,1}x(e)=\nu~_{H}D^{1,1}x(\zeta)\]
     \\\\ \textbf{AMS Classification- 34A08, 35R11}
    \end{abstract}

    \maketitle \numberwithin{equation}{section}
    \newtheorem{theorem}{Theorem}[section]
    \newtheorem{lemma}[theorem]{Lemma}
    \newtheorem{definition}[theorem]{Definition}
    \newtheorem{example}[theorem]{Example}

    \newtheorem{remark}[theorem]{Remark}
    \allowdisplaybreaks
  \[\textbf{1.Introduction}.\]
\\ \par The fractional differential equations give proofs of the more
 appropriate models for describing real world problems. Indeed, these
 problems cannot be described using classical integer order differential
 equations. In the past years the theory of fractional differential
 equations has received much attention from the authors, and has
 become an important field of investigation due to existence
 applications in engineering, biology, chemistry, economics and numerous branches of physics sciences[1,2,5,10,11].
 Fractional differential equations have a several kinds of fractional differential
 equations. One of them is the Hadamard fractional derivative innovated by
 Hadamard in 1892[3],which differs from the Riemann-Liouvill and
 Caputo type fractional derivative[10],the preceding ones in the
 sense that the kernel of the integral contains logarithmic function
 of arbitrary exponent. The properties of Hadamard Fractional integral and derivative
 can be found in[2,26].Recently, the authors studied the Hadamard-type,
 Caputo-Hadamard-type and Hilfer-Hadamard-type
 fractional derivative by using the fixed point theorems with the boundary
 value problems and give the results of existence and
 uniqueness[15-22,25].In this paper, we studied the existence and uniqueness result of
solutions for boundary value problems for Hilfer-Hadamard-Type
fractional differential equations of the form
$$_{H}D^{\alpha,\beta}x(t)+f(t,x(t))=0,~\quad\quad~~t\in
J:=(1,e],~\quad~ 1<\alpha\leq2~,\quad~~0\leq\beta\leq1~~$$
$$x(1+\epsilon)=0,~\quad\quad~_{H}D^{1,1}x(e)=\nu~_{H}D^{1,1}x(\zeta)~\quad\quad\quad\quad\quad\quad\quad\quad(1.1)$$
\\where $~_{H}D^{\alpha,\beta}~~$is the Hilfer-Hadamard fractional
derivative of order $\quad1<\alpha\leq2~~$ and
type$~\beta\in[0,1],~\quad~~0\leq\nu<1,~\quad~~\zeta\in(1,e),\quad
 0<\epsilon<1,~~\quad~_{H}D^{1,1}=t\frac{d}{dt}~$ and $~f:J\rightarrow\mathbb{R}^{+}.~$

\[\textbf{2.Preliminaries}\]In this section, we introduce some
notations and definitions of Hilfer-Hadamard-Type fractional
calculus.
\\\textbf{\ Definition 2.1.[2,11]~} (Riemann-Liouville fractional integral).
\\ The Riemann-Liouville integral of order $~\alpha~ > 0$ of a function
$~\varphi:[1,\infty)\rightarrow\mathbb{R}~$  is defined by
$$(I^{\alpha}\varphi)(t)=\frac{1}{\Gamma(\alpha)}
_{1}\int^{t}\frac{\varphi(\tau)d\tau}{(t-\tau)^{1-\alpha}}\quad,\quad(t>1),$$
Here $\Gamma(\alpha)$ is the Euler's Gamma function.
\\\textbf{\ Definition 2.2.[2,11]~} (Riemann-Liouville fractional derivative).
\\The Riemann-Liouville fractional derivative of order $~\alpha > 0~$
of a function \\$~\varphi:[1,\infty)\rightarrow\mathbb{R}~$  is
defined by \\ \par $(D^{\alpha}\varphi)(t):=(\frac{d}{dt})^{n}
(I^{n-\alpha}\varphi)(t)$
$$\quad\quad=\frac{1}{\Gamma(n-\alpha)}\frac{d^{n}}{dt^{n}}
_{1}\int^{t}\frac{\varphi(\tau)d\tau}{(t-\tau)^{\alpha-n+1}}
\quad,\quad\quad(n=[\alpha]+1;t>1),$$ Here [$\alpha$] is the integer
part of $\alpha.$
 \\\textbf{\ Definition 2.3.[2]~}(Hadamard Fractional integral).
\\~The Hadamard Fractional integral of order $~\alpha\in \mathbb{R}^{+}~$for a function
$~\varphi:[1,\infty)\rightarrow\mathbb{R}~$ is defined
as\[_{H}I^{\alpha}\varphi(t)=\frac{1}{\Gamma(\alpha)}_{1}\int^{t}(\log\frac{t}{\tau})^{\alpha-1}
\quad\frac{\varphi(\tau)}{\tau}d\tau,\quad\quad(t>1)\] where
$~\log(.)=\log_{e}(.)~$.
\\\textbf{\
Definition 2.4.[2]~}(Hadamard Fractional derivative).
\\~The Hadamard Fractional derivative of order $~\alpha~$ applied to
the function\\ $~\varphi:[1,\infty)\rightarrow\mathbb{R}~$ is
defined
as\[_{H}D^{\alpha}\varphi(t)=\delta^{n}(_{H}I^{n-\alpha}\varphi(t)),\quad
n-1<\alpha<n,\quad n=[\alpha]+1,\]
where$\quad~\delta^{n}=(t\frac{d}{dt})^{n}\quad~$and$~[\alpha]~$denotes
the integer part of the real number$~\alpha.~$
\\\textbf{\ Definition 2.5.[3,13]~}(Caputo-Hadamard Fractional derivative).
\\~The Caputo-Hadamard Fractional derivative of
order$~\alpha~$applied to the function $~\varphi\in
AC_{\delta}^{n}[a,b]~$is defined
as\[_{HC}D^{\alpha}\varphi(t)=(_{H}I^{n-\alpha}\delta^{n}\varphi)(t)\]
where\\$~n=[\alpha]+1,~$and$~\varphi\in
AC_{\delta}^{n}[a,b]=\bigg\{\varphi:[a,b]\rightarrow
\mathbb{C}:\delta^{(n-1)}\varphi\in
AC[a,b],\delta=t\frac{d}{dt}\bigg\}~$
\\\textbf{\
Definition 2.6.[5,21]~}(Hilfer Fractional derivative).
\\Let$~~n-1<\alpha<n,~~0\leq\beta\leq 1,~~\varphi\in
L^{1}(a,b).~$The Hilfer Fractional derivative $D^{\alpha,\beta}$ of
order$~\alpha~$and type $~\beta~$of$~\varphi~$is defined
as\[~(D^{\alpha,\beta}\varphi)(t)=\big(I^{\beta(n-\alpha)}(\frac{d}{dt})^{n}~I^{(n-\alpha)(1-\beta)}\varphi\big)(t)~\]
\[=\big(I^{\beta(n-\alpha)}(\frac{d}{dt})^{n}I^{n-\gamma}\varphi\big)(t);~\quad\gamma=\alpha+n\beta-\alpha\beta.\]
\[=\big(I^{\beta(n-\alpha)}D^{\gamma}\varphi\big)(t),\]
\\Where $I^{(.)}$ and  $~D^{(.)}~$is the Riemann-Liouvill fractional integral and
derivative defined by (2.1) and (2.2), respectively.\\In particular,
if $\quad0<\alpha<1,$ then
$$~(D^{\alpha,\beta}\varphi)(t)=\big(I^{\beta(1-\alpha)}\frac{d}{dt}~I^{(1-\alpha)(1-\beta)}\varphi\big)(t)~$$
\[\quad\quad=\big(I^{\beta(1-\alpha)}\frac{d}{dt}I^{1-\gamma}\varphi\big)(t);\quad\gamma=\alpha+\beta-\alpha\beta.~\]
\[=\big(I^{\beta(1-\alpha)}D^{\gamma}\varphi\big)(t).\]
\\ \textbf {\ Properties 2.7.[21,22].}\\ Let$~~0<\alpha<1,\quad0\leq\beta\leq 1,
\quad\gamma=\alpha+\beta-\alpha\beta,$ and $\varphi\in L^{1}(a,b).~$
If $D^{\gamma}\varphi$ exists and in $L^{1}(a,b),$ then
\[I_{a+}^{\alpha}~(D_{a+}^{\alpha,\beta}\varphi)(t)=I_{a+}^{\gamma}~
 (D_{a+}^{\gamma}\varphi)(t)=
 \varphi(t)-\frac{(I_{a+}^{1-\gamma}\varphi)(a)}{\Gamma(\gamma)}(t-a)^{\gamma-1}\]
\\\textbf{\ Definition 2.8.[21,22]}(Hilfer-Hadamard Fractional derivative).
\\~Let$~~0<\alpha<1,~~0\leq\beta\leq 1,~~\varphi\in
L^{1}(a,b).~$The Hilfer-Hadamard Fractional derivative
$_{H}D^{\alpha,\beta}$ of order$~\alpha~$and type
$~\beta~$of$~\varphi~$is defined
as\[~(_{H}D^{\alpha,\beta}\varphi)(t)=\big(_{H}I^{\beta(1-\alpha)}\delta~_{H}I^{(1-\alpha)(1-\beta)}\varphi\big)(t)~\]
\[=\big(_{H}I^{\beta(1-\alpha)}\delta~_{H}I^{1-\gamma}\varphi\big)(t);\quad\gamma=\alpha+\beta-\alpha\beta.\]
\[=\big(_{H}I^{\beta(1-\alpha)}_{H}D^{\gamma}\varphi\big)(t).\]
Where $_{H}I^{(.)}$ and $~_{H}D^{(.)}~$is the Hadamard fractional
integral and derivative defined by (2.3) and (2.4), respectively.
\\ \textbf {\ Theorem 2.9.[2,3].}
 \\~Let$~\Re(\alpha)>0,\quad~n=[\Re(\alpha)]+1~$and$~0<a<b<\infty.~$ if $~\varphi\in L^{1}(a,b)~$ and
 $~(_{H}I_{a+}^{n-\alpha}\varphi)(t)\in AC_{\delta}^{n}[a,b],~$ then
 \[(_{H}I_{a+}^{\alpha}~_{H}D_{a+}^{\alpha}\varphi)(t)=
 \varphi(t)-\sum_{j=0}^{n-1}\frac{(\delta^{(n-j-1)}(_{H}I_{a+}^{n-\alpha}\varphi))(a)}{\Gamma(\alpha-j)}(\log\frac{t}{a})^{\alpha-j-1}\]
\\\textbf{\ Theorem 2.10.[3,13]~}\\Let  $\varphi(t)\in AC_{\delta}^{n}[a,b]$ or $\varphi(t)\in
C_{\delta}^{n}[a,b],$ and $~~\alpha\in\mathbb{C},$
then\[(_{H}I_{a+}^{\alpha}~_{HC}D_{a+}^{\alpha}\varphi)(t)=
 \varphi(t)-\sum_{K=0}^{n-1}\frac{\delta^{K}\varphi(a)}{\Gamma(K+1)}(\log\frac{t}{a})^{K}\]
 \\\textbf{\ Theorem 2.11.[23,24]~}(Leray- Schauder alternative).
\\Let $E$ be a Banach space. Suppose that $T:E\rightarrow E$ is
completely continuous operator and the set $V = \big\{v\in E\mid v
 =\lambda  Tv, 0 < \lambda < 1\big\}$ is bounded. Then $T$ has a fixed point
in $E$.
\\\textbf{\ Theorem 2.12}[23,24]
\\Let $E$ be a Banach space and  $V$ is an open bounded
subset of $E$ with $0\in V$.\\Suppose that $\Psi :
\overline{V}\rightarrow E$ be a completely continuous operator such
that \\$\parallel\Psi v\parallel\leq\parallel v\parallel,  \forall
v\in \partial V$. Then $\Psi$ has a fixed point in $\overline{V}.$
 \\\textbf{\ Lemma 2.13.[25]} \\For $~~~1<\alpha\leq2~~~~~$and $~~~~\varphi\in
C([1,e]~,~\mathbb{R})~$ the problem for Capoto-Hadamard-type,
$$~~~~_{HC}D^{\alpha}x(t)+\varphi(t)=0,~~\quad\quad t\in [1,e]~\quad~
1<\alpha\leq2,~$$~ $$x(1)=0,\quad\quad_{HC}D x(e)=\nu~_{HC}D
x(\zeta)~$$
\\has a unique solution it giving in the formulae
\\\[x(t)=-\frac{1}{\Gamma(\alpha)}_{1}\int^{t}(\log\frac{t}{\tau})^{\alpha-1}\frac{\varphi(\tau)}{\tau}d\tau~\quad\quad\quad\quad\quad\quad\quad\]
\[\quad\quad+\frac{\log t }{1-\nu}
\bigg[\frac{1}{\Gamma(\alpha-1)}_{1}\int^{e}(\log\frac{e}{\tau})^{\alpha-2}\frac{\varphi(\tau)}{\tau}d\]\[\quad\quad\quad-
\frac{\nu}{\Gamma(\alpha-1)}_{1}\int^{\zeta}(\log\frac{\zeta}{\tau})^{\alpha-2}\frac{\varphi(\tau)}{\tau}d\tau\bigg]\]
\\\[\textbf{3.Main Results }.\]
\\\textbf{\ Definition 3.1}(Hilfer-Hadamard Fractional derivative).
\\~Let$~~n-1<\alpha<n,~~0\leq\beta\leq 1,~~\varphi\in
L^{1}(a,b).~$The Hilfer-Hadamard Fractional derivative
$_{H}D^{\alpha,\beta}$ of order$~\alpha~$and type
$~\beta~$of$~\varphi~$is defined
as\[~(_{H}D^{\alpha,\beta}\varphi)(t)=\big(_{H}I^{\beta(n-\alpha)}(\delta)^{n}~_{H}I^{(n-\alpha)(1-\beta)}\varphi\big)(t)~\]
\[=\big(_{H}I^{\beta(n-\alpha)}(\delta)^{n}~_{H}I^{n-\gamma}\varphi\big)(t);\quad\gamma=\alpha+n\beta-\alpha\beta.\]
\[=\big(I^{\beta(n-\alpha)}_{H}D^{\gamma}\varphi\big)(t),\]
Where $_{H}I^{(.)}$ and $~_{H}D^{(.)}~$is the Hadamard fractional
integral and derivative defined by (2.3) and (2.4), respectively.
\\ \textbf{\ Lemma 3.2.}
 \\~Let$~\Re(\alpha)>0,\quad0\leq\beta\leq1,\quad\gamma=\alpha+n\beta-\alpha\beta,
 \quad n-1<\gamma\leq n,\quad~n=[\Re(\alpha)]+1~$ and $~0<a<b<\infty.~$ if $~\varphi\in L^{1}(a,b)~$ and
 $~(_{H}I_{a+}^{n-\gamma}\varphi)(t)\in AC_{\delta}^{n}[a,b],
\quad ~~$ then
 \[_{H}I_{a+}^{\alpha}~(_{H}D_{a+}^{\alpha,\beta}\varphi)(t)=_{H}I_{a+}^{\gamma}~(_{H}D_{a+}^{\gamma}\varphi)(t)=
 \varphi(t)-\sum_{j=0}^{n-1}\frac{(\delta^{(n-j-1)}(_{H}I_{a+}^{n-\gamma}\varphi))(a)}{\Gamma(\gamma-j)}(\log\frac{t}{a})^{\gamma-j-1}\]
 From this Lemma, we notice that if $\beta=0$ the formulae reduces to the formulae in
 the theorem 2.9, and if the $\beta=1$ the formulae reduces to the formulae in
 the theorem 2.10.
\\ \textbf{\ Proof.}  We have\\
\\$_{H}I_{a+}^{\alpha}~(_{H}D_{a+}^{\alpha,\beta}\varphi)(t)
=_{H}I_{a+}^{\gamma}~(_{H}D_{a+}^{\gamma}\varphi)(t)$
\begin{align*}
&=\frac{1}{\Gamma(\gamma)}_{a}\int^{t}(\log\frac{t}{\tau})^{\gamma-1}
(_{H}D_{a+}^{\gamma}\varphi(\tau))\frac{d\tau}{\tau}\\
&=\frac{d}{dt}
\biggr\{\frac{1}{\Gamma(\gamma+1)}_{a}\int^{t}(\log\frac{t}{\tau})^{\gamma}
(_{H}D_{a+}^{\gamma}\varphi(\tau))\frac{d\tau}{\tau}\biggr\}
\end{align*}
\par On the hand, repeatedly integrating by parts and then using
\\$_{H}I_{a+}^{p} . _{H}I_{a+}^{q}=_{H}I_{a+}^{q} .
_{H}I_{a+}^{p}=_{H}I_{a+}^{p+q},$ we obtain
\begin{align*}
&\frac{1}{\Gamma(\gamma+1)}_{a}\int^{t}(\log\frac{t}{\tau})^{\gamma}
(_{H}D_{a+}^{\gamma}\varphi(\tau))\frac{d\tau}{\tau}\\&\quad\quad=
\frac{1}{\Gamma(\gamma+1)}_{a}\int^{t}(\log\frac{t}{\tau})^{\gamma}
\delta^{n}(_{H}I_{a+}^{n-\gamma}\varphi(\tau))\frac{d\tau}{\tau}
\\&\quad\quad=-\frac{1}{\Gamma(\gamma+1)}(\log\frac{t}{a})^{\gamma}\delta^{n-1}(_{H}I_{a+}^{n-\gamma}\varphi(a))+
\frac{1}{\Gamma(\gamma)}_{a}\int^{t}(\log\frac{t}{\tau})^{\gamma-1}
\delta^{n-1}(_{H}I_{a+}^{n-\gamma}\varphi(\tau))\frac{d\tau}{\tau}\\
&\quad\quad=-\frac{1}{\Gamma(\gamma+1)}(\log\frac{t}{a})^{\gamma}\delta^{n-1}(_{H}I_{a+}^{n-\gamma}\varphi(a))
-\frac{1}{\Gamma(\gamma)}(\log\frac{t}{a})^{\gamma-1}\delta^{n-2}(_{H}I_{a+}^{n-\gamma}\varphi(a))
\\&\quad\quad\quad+\frac{1}{\Gamma(\gamma-1)}_{a}\int^{t}(\log\frac{t}{\tau})^{\gamma-2}
\delta^{n-2}(_{H}I_{a+}^{n-\gamma}\varphi(\tau))\frac{d\tau}{\tau}\\
&\quad\quad=-\frac{1}{\Gamma(\gamma+1)}(\log\frac{t}{a})^{\gamma}\delta^{n-1}(_{H}I_{a+}^{n-\gamma}\varphi(a))
-\frac{1}{\Gamma(\gamma)}(\log\frac{t}{a})^{\gamma-1}\delta^{n-2}(_{H}I_{a+}^{n-\gamma}\varphi(a))\\
&\quad\quad\quad-\frac{1}{\Gamma(\gamma-1)}(\log\frac{t}{a})^{\gamma-2}\delta^{n-3}
(_{H}I_{a+}^{n-\gamma}\varphi(a))+\frac{1}{\Gamma(\gamma-2)}_{a}\int^{t}(\log\frac{t}{\tau})^{\gamma-3}
\delta^{n-3}(_{H}I_{a+}^{n-\gamma}\varphi(\tau))\frac{d\tau}{\tau}\\&\quad\quad=
\frac{1}{\Gamma(\gamma-2)}_{a}\int^{t}(\log\frac{t}{\tau})^{\gamma-3}
\delta^{n-3}(_{H}I_{a+}^{n-\gamma}\varphi(\tau))\frac{d\tau}{\tau}-\sum_{j=1}^{3}
\frac{\delta^{n-j}(_{H}I_{a+}^{n-\gamma}\varphi(a))}{\Gamma(2+\gamma-j)}(\log\frac{t}{a})^{\gamma-j+1}
\\&\quad\quad=......\\&\quad\quad=
\frac{1}{\Gamma(\gamma-n+1)}_{a}\int^{t}(\log\frac{t}{\tau})^{\gamma-n}
(_{H}I_{a+}^{n-\gamma}\varphi(\tau))\frac{d\tau}{\tau}-\sum_{j=1}^{n}
\frac{\delta^{n-j}(_{H}I_{a+}^{n-\gamma}\varphi(a))}{\Gamma(2+\gamma-j)}(\log\frac{t}{a})^{\gamma-j+1}
\end{align*}
Therefore,\\
\\$_{H}I_{a+}^{\alpha}~(_{H}D_{a+}^{\alpha,\beta}\varphi)(t)
=_{H}I_{a+}^{\gamma}~(_{H}D_{a+}^{\gamma}\varphi)(t)$
\begin{align*}
&=\frac{d}{dt} \biggr\{(_{H}I_{a+}^{1}\varphi)(t)-\sum_{j=1}^{n}
\frac{\delta^{n-j}(_{H}I_{a+}^{n-\gamma}\varphi(a))}{\Gamma(2+\gamma-j)}(\log\frac{t}{a})^{\gamma-j+1}\biggr\}\\
\\&=\varphi(t)-\sum_{j=0}^{n-1}
\frac{\delta^{n-j-1}(_{H}I_{a+}^{n-\gamma}\varphi(a))}{\Gamma(\gamma-j)}
(\log\frac{t}{a})^{\gamma-j-1}\quad\quad\quad\quad\quad\quad\quad\quad\quad\quad\quad\quad\Box
\end{align*}
 \\\textbf{\ Lemma 3.3.}
 \\For $1<\alpha\leq2~~$,$~~0\leq\beta\leq1~~$ and $~~\varphi\in
C([1,e],\mathbb{R}),$\\$~~~\gamma=\alpha+2\beta-\alpha\beta~~$,$~\gamma\in(1,2]$
\\ the problem
\\$_{H}D^{\alpha,\beta}x(t)+\varphi(t)=0,~~~~~~~~~~~~$ $t\in J~~,$ $1<\alpha\leq2~~$,$~~0\leq\beta\leq1~~$
\\$x(1+\epsilon)=0,\quad\quad\quad~$$_{H}D^{1,1}x(e)=\nu~_{H}D^{1,1}x(\zeta)\quad\quad\quad\quad\quad \quad\quad\quad\quad\quad\quad\quad\quad(3.1)$
\\has a unique solution it giving in the formulae
\\\[x(t)=-\frac{1}{\Gamma(\alpha)}_{1}\int^{t}(\log\frac{t}{\tau})^{\alpha-1}\frac{\varphi(\tau)}{\tau}d\tau\quad\quad\quad
\quad\quad\quad\quad\quad\quad\quad\quad\quad\quad\quad~\]\[
 +(\frac{\log t}{\log(1+\epsilon)})^{\gamma-1}~\frac{1}{\Gamma(\alpha)}_{1}\int^{1+\epsilon}(\log\frac{1+\epsilon}{\tau})^{\alpha-1}
 \frac{\varphi(\tau)}{\tau}d\tau~~~\]
 \[\quad\quad\quad~+[\frac{1}{\log(1+\epsilon)}-\frac{1}{\log t}]\frac{(\log
t)^{\gamma-1}}{\sum_{i=0}^1\eta_{i}}\bigg[\frac{1}{\Gamma(\alpha-1)}_{1}\int^{e}(\log\frac{e}{\tau})^{\alpha-2}\frac{\varphi(\tau)}{\tau}d\tau\]
\[~~~~~~~~-\frac{\nu}{\Gamma(\alpha-1)}_{1}\int^{\zeta}(\log\frac{\zeta}{\tau})^{\alpha-2}\frac{\varphi(\tau)}{\tau}d\tau
+\varepsilon[1-\nu(\log \zeta)^{\gamma-2}]\bigg]\]
\\Where
\begin{align*}
&\sum_{i=0}^1\eta_{i}=\sum_{i=0}^1(-1)^{i}(\log(1+\epsilon))^{i-1}(\gamma-i-1)[1-\nu(\log\zeta)^{\gamma-i-2}],\quad
with\quad\sum_{i=0}^1\eta_{i}\neq0
\\&\varepsilon=(1-\gamma)(\log(1+\epsilon))^{1-\gamma}~_{H}I^{\alpha}\varphi(1+\epsilon)~~~~~
\end{align*}
\\ \textbf{\ Proof.} In the view of the Lemma$(3.2)$, the solution of
the Hilfer-Hadamard differential equation $(3.1)$ can be written as
\[\quad\quad\quad\quad x(t)=-~_{H}I^{\alpha}\varphi(t)+c_{0}(\log t)^{\gamma-1}+c_{1}(\log t)^{\gamma-2}~
\quad\quad\quad\quad\quad\quad\quad\quad\quad\quad\quad~(3.2)\]
and\[\quad\quad\quad _{H}D^{1,1}
x(t)=-~_{H}I^{\alpha-1}\varphi(t)+(\gamma-1)c_{0}(\log
t)^{\gamma-2}+(\gamma-2) c_{1}(\log
t)^{\gamma-3}~\quad\quad\quad~(3.3)\] The boundary condition
$~~x(1+\epsilon)=0~~$gives
\[\quad\quad c_{0}=(\log
(1+\epsilon))^{1-\gamma}~_{H}I^{\alpha}\varphi(1+\epsilon)-\frac{c_{1}}{\log
(1+\epsilon)}~\quad\quad\quad\quad\quad\quad\quad\quad\quad\quad\quad\quad\quad\quad~(3.4)\]
In view of the boundary condition
$~~_{H}D^{1,1}x(e)=\nu~_{H}D^{1,1}x(\zeta),~~$ and by $~(3.3),~$ and
$~(3.4),~$  we have
\[c_{1}=\frac{1}{\sum_{i=0}^1\eta_{i}}\bigg[-_{H}I^{\alpha-1}\varphi(e)+\nu~_{H}I^{\alpha-1}\varphi(\zeta)+\varepsilon[1-\nu(\log
\zeta)^{\gamma-2}]\bigg]\]
\\Where
\begin{align*}
&\eta_{i}=(-1)^{i}(\log(1+\epsilon))^{i-1}(\gamma-i-1)[1-\nu(\log\zeta)^{\gamma-i-2}],
\\&\varepsilon=(1-\gamma)(\log(1+\epsilon))^{1-\gamma}~_{H}I^{\alpha}\varphi(1+\epsilon)~~~~~
\end{align*}
 Substituting the value of $~c_{1}~$in$~(3.4)~$ we have
\begin{align*}
~c_{0}&=(\log(1+\epsilon))^{1-\gamma}~_{H}I^{\alpha}\varphi(1+\epsilon)\\&\quad\quad\quad-\frac{1}{\sum_{i=0}^
1\eta_{i}\log(1+\epsilon)}\bigg[-_{H}I^{\alpha-1}\varphi(e)\\&\quad\quad\quad\quad\quad\quad
+\nu~_{H}I^{\alpha-1}\varphi(\zeta)+\varepsilon[1-\nu(\log
\zeta)^{\gamma-2}]\bigg]~
\end{align*}
Now substituting the values of $~c_{0}~$ and $~c_{1}~$in $~(3.2)~$we
obtain the solution of the problem(3.1).
\[\textbf{Results of Existence}.\]
Suppose that
$$\quad\quad\quad\quad\quad\quad\quad\quad\quad\quad\quad~K=C([1,e],\mathbb{R})\quad\quad\quad\quad\quad\quad
\quad\quad\quad\quad\quad\quad\quad\quad\quad\quad\quad\quad(3.5)$$
is a Banach space of all continuous functions from $[1,e]$ into
$\mathbb{R}$ talented with the norm $\parallel x\parallel=\sup_{t\in
J}\mid x(t)\mid.$
\\\par From the Lemma3.1,we getting an operator $~\rho:K\rightarrow K$
defined as
\begin{align*}
\quad\quad\quad(\rho
x)(t)&=-\frac{1}{\Gamma(\alpha)}_{1}\int^{t}(\log\frac{t}{\tau})^{\alpha-1}\frac{\varphi(\tau)}{\tau}d\tau
\\&\quad\quad\quad +(\frac{\log t}{\log(1+\epsilon)})^{\gamma-1}~\frac{1}{\Gamma(\alpha)}
_{1}\int^{1+\epsilon}(\log\frac{1+\epsilon}{\tau})^{\alpha-1}\frac{\varphi(\tau)}{\tau}d\tau\quad\quad\quad\quad(3.6)
\\&\quad+[\frac{1}{\log(1+\epsilon)}-\frac{1}{\log t}]\frac{(\log
t)^{\gamma-1}}{\sum_{i=0}^1\eta_{i}}\bigg[\frac{1}{\Gamma(\alpha-1)}_{1}\int^{e}(\log\frac{e}{\tau})^{\alpha-2}\frac{\varphi(\tau)}{\tau}d\tau
\\&\quad\quad\quad-\frac{\nu}{\Gamma(\alpha-1)}_{1}\int^{\zeta}(\log\frac{\zeta}{\tau})^{\alpha-2}\frac{\varphi(\tau)}{\tau}d\tau
+\varepsilon[1-\nu(\log \zeta)^{\gamma-2}]\bigg],\quad t\in J~
\end{align*}
 It must be noticed that
the problem $(1.1)$ has solutions if and only if the operator $\rho$
has fixed points.The result of existence and uniqueness is based on
the Banach Principle of contraction.
\\\textbf{\ Theorem 3.4}  suppose that there exists a constant
$~C>0~$such that \\$~~\mid f(t,x(t))-f(t,y(t))\mid\leq C\mid
x-y\mid,~~\forall t\in J,~~C>0,x,y\in \mathbb{R}.~~~~~~~~~$
\\If$~\Phi$ satisfied the condition  $~~C\Phi<1,~$ where
\begin{align*}
\quad\quad\quad\Phi&=\bigg\{\frac{[1+(\log(1+\epsilon))^{1-\gamma+\alpha}]}{\Gamma(\alpha+1)}\\&\quad\quad\quad\quad+
\frac{[1-(\log(1+\epsilon))]}{(\log(1+\epsilon))\Gamma(\alpha)\sum_{i=0}^1\eta_{i}}\bigg[1+\nu(\log
\zeta)^{\alpha-1}\\&\quad\quad\quad\quad\quad\quad\quad\quad+\frac{(1-\gamma)(\log(1+\epsilon))^{1-\gamma+\alpha}}{\alpha}[1-\nu(\log
\zeta)^{\gamma-2}]\bigg]\bigg\}\quad\quad\quad\quad(3.7)
\end{align*}
Then the problem $(1.1)$ has a unique solution.
\\\textbf{\ Proof.} We put$~~sup_{t\in J}\mid f(\tau,0)\mid=P<\infty~$ and choose
$~r\geq\frac{\Phi P}{1-\Phi C}.$
\par Now,assume that $B_{r}=\{x\in K:\parallel x\parallel\leq r\},$
then we show that $\rho B_{r}\subset B_{r}.$ \\For any $x\in B_{r},$
we have
\begin{align*}
&\parallel(\rho x)(t)\parallel\\&\quad=sup_{t\in J}\bigg\{\Biggl|
-\frac{1}{\Gamma(\alpha)}_{1}\int^{t}(\log\frac{t}{\tau})^{\alpha-1}f(\tau,x(\tau))\frac{d\tau}{\tau}~
\\&\quad\quad\quad\quad+(\frac{\log t}{\log(1+\epsilon)})^{\gamma-1}~\frac{1}{\Gamma(\alpha)}_{1}\int^{1+\epsilon}
(\log\frac{1+\epsilon}{\tau})^{\alpha-1}f(\tau,x(\tau))\frac{d\tau}{\tau}
\\&\quad\quad+\big[\frac{1}{\log(1+\epsilon)}-\frac{1}{\log t}\big]\frac{(\log
t)^{\gamma-1}}{\sum_{i=0}^1\eta_{i}}\bigg[\frac{1}
{\Gamma(\alpha-1)}_{1}\int^{e}(\log\frac{e}{\tau})^{\alpha-2}f(\tau,x(\tau))\frac{d\tau}{\tau}
\\&\quad\quad\quad\quad-\frac{\nu}{\Gamma(\alpha-1)}_{1}\int^{\zeta}
(\log\frac{\zeta}{\tau})^{\alpha-2}f(\tau,x(\tau))\frac{d\tau}{\tau}
\\&\quad\quad+\frac{(1-\gamma)(\log(1+\epsilon))^{1-\gamma}[1-\nu(\log
\zeta)^{\gamma-2}]}{\Gamma(\alpha)}_{1}\int^{1+\epsilon}
(\log\frac{1+\epsilon}{\tau})^{\alpha-1}f(\tau,x(\tau))\frac{d\tau}{\tau}\bigg]\Biggl|\bigg\}
\\&\quad\leq\frac{1}{\Gamma(\alpha)}_{1}\int^{t}(\log\frac{t}{\tau})^{\alpha-1}
\bigg(\mid f(\tau,x(\tau))-f(\tau,0)\mid+\mid
f(\tau,0)\mid\bigg)\frac{d\tau}{\tau}~\\&\quad\quad\quad+(\frac{\log
t}{\log(1+\epsilon)})^{\gamma-1}~\frac{1}{\Gamma(\alpha)}_{1}\int^{1+\epsilon}(\log\frac{1+\epsilon}{\tau})^{\alpha-1}
\bigg(\mid f(\tau,x(\tau))-f(\tau,0)\mid+\mid
f(\tau,0)\mid\bigg)\frac{d\tau}{\tau}
\\&\quad+\big[\frac{1}{\log(1+\epsilon)}-\frac{1}{\log t}\big]\frac{(\log
t)^{\gamma-1}}{\sum_{i=0}^1\eta_{i}}\bigg[\frac{1}{\Gamma(\alpha-1)}_{1}\int^{e}(\log\frac{e}{\tau})^{\alpha-2}\bigg(\mid
f(\tau,x(\tau))-f(\tau,0)\mid+\mid
f(\tau,0)\mid\bigg)\frac{d\tau}{\tau}
\\&\quad\quad\quad\quad+\frac{\nu}{\Gamma(\alpha-1)}_{1}\int^{\zeta}(\log\frac{\zeta}{\tau})^{\alpha-2}\bigg(\mid
f(\tau,x(\tau))-f(\tau,0)\mid+\mid
f(\tau,0)\mid\bigg)\frac{d\tau}{\tau}
\\&\quad+\frac{(1-\gamma)(\log(1+\epsilon))^{1-\gamma}[1-\nu(\log
\zeta)^{\gamma-2}]}{\Gamma(\alpha)}
_{1}\int^{1+\epsilon}(\log\frac{1+\epsilon}{\tau})^{\alpha-1}\bigg(\mid
f(\tau,x(\tau))-f(\tau,0)\mid+\mid
f(\tau,0)\mid\bigg)\frac{d\tau}{\tau}\bigg]
\\&\leq(Cr+P)\bigg\{\frac{1}{\Gamma(\alpha+1)}+\frac{(\log(1+\epsilon))^{1-\gamma+\alpha}}{\Gamma(\alpha+1)}
\\&\quad\quad\quad\quad\quad\quad\quad+\frac{[1-
\log(1+\epsilon)]}{\log(1+\epsilon)\sum_{i=0}^1\eta_{i}\Gamma(\alpha)}\bigg[1+\nu(\log\zeta)^{\alpha-1}\\&\quad\quad\quad\quad\quad\quad
\quad+\frac{(1-\gamma)(\log(1+\epsilon))^{1-\gamma+\alpha}
[1-\nu(\log\zeta)^{\gamma-2}]}{\alpha}\bigg]\bigg\}
\\&\leq(Cr+P)\bigg\{\frac{[1+(\log(1+\epsilon))^{1-\gamma+\alpha}]}{\Gamma(\alpha+1)}\\&\quad\quad\quad\quad\quad\quad+
\frac{[1-(\log(1+\epsilon))]}{(\log(1+\epsilon))\Gamma(\alpha)\sum_{i=0}^1\eta_{i}}\bigg[1+\nu(\log
\zeta)^{\alpha-1}\\&\quad\quad\quad\quad\quad\quad\quad\quad+\frac{(1-\gamma)(\log(1+\epsilon))^{1-\gamma+\alpha}}{\alpha}[1-\nu(\log
\zeta)^{\gamma-2}]\bigg]\bigg\}
\\&\leq(Cr+P)\Phi\leq
r\quad\quad\quad\quad\quad\quad\quad\quad\quad\quad\quad\quad\quad\quad\quad\quad\quad\quad
\quad\quad\quad\quad\quad\quad\quad\quad\quad\quad\quad\quad\quad\quad(3.8)
\end{align*}
Thus we shown $\rho B_{r}\subset B_{r}.$
\\Now,For $x,y\in K$ and $\forall t\in J,$ we have
\begin{align*}
&\mid(\rho x)(t)-(\rho
y)(t)\mid\\&\quad=\Biggl|-\frac{1}{\Gamma(\alpha)}_{1}\int^{t}(\log\frac{t}{\tau})^{\alpha-1}
\bigg(f(\tau,x(\tau))-f(\tau,y(\tau))\bigg)\frac{d\tau}{\tau}~
\\&\quad\quad\quad+(\frac{\log t}{\log(1+\epsilon)})^{\gamma-1}~\frac{1}{\Gamma(\alpha)}_{1}\int^{1+\epsilon}
(\log\frac{1+\epsilon}{\tau})^{\alpha-1}\bigg(f(\tau,x(\tau))-f(\tau,y(\tau))\bigg)\frac{d\tau}{\tau}
\\&\quad\quad\quad+[\frac{1}{\log(1+\epsilon)}-\frac{1}{\log t}]\frac{(\log
t)^{\gamma-1}}{\sum_{i=0}^1\eta_{i}}\bigg[\frac{1}{\Gamma(\alpha-1)}_{1}\int^{e}
(\log\frac{e}{\tau})^{\alpha-2}\bigg(f(\tau,x(\tau))-f(\tau,y(\tau))\bigg)\frac{d\tau}{\tau}
\\&\quad\quad\quad\quad\quad\quad\quad-\frac{\nu}{\Gamma(\alpha-1)}_{1}\int^{\zeta}
(\log\frac{\zeta}{\tau})^{\alpha-2}\bigg(f(\tau,x(\tau))-f(\tau,y(\tau))\bigg)\frac{d\tau}{\tau}
\\&\quad\quad\quad+\frac{(1-\gamma)(\log(1+\epsilon))^{1-\gamma}[1-\nu(\log
\zeta)^{\gamma-2}]}{\Gamma(\alpha)}_{1}\int^{1+\epsilon}
(\log\frac{1+\epsilon}{\tau})^{\alpha-1}\bigg(f(\tau,x(\tau))-f(\tau,y(\tau))\bigg)\frac{d\tau}{\tau}\bigg]\Biggl|
\\&\quad\leq\frac{1}{\Gamma(\alpha)}_{1}\int^{t}(\log\frac{t}{\tau})^{\alpha-1}\mid f(\tau,x(\tau))-f(\tau,y(\tau))\mid\frac{d\tau}{\tau}~
\\&\quad\quad\quad+(\frac{\log t}{\log(1+\epsilon)})^{\gamma-1}~\frac{1}{\Gamma(\alpha)}_{1}\int^{1+\epsilon}
(\log\frac{1+\epsilon}{\tau})^{\alpha-1}\mid
f(\tau,x(\tau))-f(\tau,y(\tau))\mid\frac{d\tau}{\tau}
\\&\quad\quad\quad+[\frac{1}{\log(1+\epsilon)}-\frac{1}{\log t}]\frac{(\log
t)^{\gamma-1}}{\sum_{i=0}^1\eta_{i}}\bigg[\frac{1}{\Gamma(\alpha-1)}_{1}\int^{e}(\log\frac{e}{\tau})^{\alpha-2}\mid
f(\tau,x(\tau))-f(\tau,y(\tau))\mid\frac{d\tau}{\tau}
\\&\quad\quad\quad\quad\quad\quad\quad+\frac{\nu}{\Gamma(\alpha-1)}_{1}\int^{\zeta}
(\log\frac{\zeta}{\tau})^{\alpha-2}\mid
f(\tau,x(\tau))-f(\tau,y(\tau))\mid\frac{d\tau}{\tau}
\\&\quad\quad\quad+\frac{(1-\gamma)(\log(1+\epsilon))^{1-\gamma}[1-\nu(\log
\zeta)^{\gamma-2}]}{\Gamma(\alpha)}_{1}\int^{1+\epsilon}(\log\frac{1+\epsilon}{\tau})^{\alpha-1}\mid
f(\tau,x(\tau))-f(\tau,y(\tau))\mid\frac{d\tau}{\tau}\bigg]
\\&\quad\leq C\parallel x-y \parallel\bigg\{\frac{[1+(\log(1+\epsilon))^{1-\gamma+\alpha}]}{\Gamma(\alpha+1)}\\&\quad\quad\quad\quad\quad\quad\quad\quad+
\frac{[1-(\log(1+\epsilon))]}{(\log(1+\epsilon))\Gamma(\alpha)\sum_{i=0}^1\eta_{i}}\bigg[1+\nu(\log
\zeta)^{\alpha-1}\\&\quad\quad\quad\quad\quad\quad\quad\quad\quad\quad\quad\quad+\frac{(1-\gamma)(\log(1+\epsilon))^{1-\gamma+\alpha}}{\alpha}[1-\nu(\log
\zeta)^{\gamma-2}]\bigg]\bigg\}
\\&\quad\leq C\Phi\parallel x-y
\parallel\quad\quad\quad\quad\quad\quad\quad\quad\quad\quad\quad\quad\quad\quad
\quad\quad\quad\quad\quad\quad\quad\quad\quad\quad\quad\quad\quad\quad\quad\quad\quad\quad(3.9)
\end{align*}
Therefore it shown that $\parallel (\rho x)(t)-(\rho y)(t)
\parallel\leq C \Phi\parallel x-y \parallel,$ where $C \Phi<1.$
\\Hence $\rho$ is a contraction. Thus by the mapping of contraction
principle the problem $(1.1)$ has a uniqueness solution.
\\ \textbf{\ Theorem 3.5}  suppose that there exists a constant $~C_{1}>0~$such
that \\$~~\mid f(t,x(t))\mid\leq C_{1},$for each $t\in J, x\in
\mathbb{R}.$ Then the problem $(1.1)$ has at least one solution.
\\\textbf{\ Proof.} The proof of this theorem will be given in
several steps, firstly, we will show that the operator $\rho$ is
completely continuous for this, in the view of the continuity of
$f$,we note that the operator $\rho$ is continuous.
\\Now, Assume that $\rho \subset \rho$ be a bounded set. \\By the
supposition that $~~\mid f(t,x(t))\mid\leq C_{1}$,for each $t\in J,
x\in \mathbb{R},$ we get
\begin{align*}
&\mid (\rho x)(t)\mid\\& \quad=\Biggl|
-\frac{1}{\Gamma(\alpha)}_{1}\int^{t}(\log\frac{t}{\tau})^{\alpha-1}f(\tau,x(\tau))\frac{d\tau}{\tau}~
\\&\quad\quad\quad+(\frac{\log t}{\log(1+\epsilon)})^{\gamma-1}~\frac{1}{\Gamma(\alpha)}_{1}\int^{1+\epsilon}
(\log\frac{1+\epsilon}{\tau})^{\alpha-1}f(\tau,x(\tau))\frac{d\tau}{\tau}
\\&\quad\quad\quad+[\frac{1}{\log(1+\epsilon)}-\frac{1}{\log t}]\frac{(\log
t)^{\gamma-1}}{\sum_{i=0}^1\eta_{i}}\bigg[\frac{1}{\Gamma(\alpha-1)}
_{1}\int^{e}(\log\frac{e}{\tau})^{\alpha-2}f(\tau,x(\tau))\frac{d\tau}{\tau}
\\&\quad\quad\quad\quad\quad\quad-\frac{\nu}{\Gamma(\alpha-1)}
_{1}\int^{\zeta}(\log\frac{\zeta}{\tau})^{\alpha-2}f(\tau,x(\tau))\frac{d\tau}{\tau}
\\&\quad\quad\quad+\frac{(1-\gamma)(\log(1+\epsilon))^{1-\gamma}[1-\nu(\log
\zeta)^{\gamma-2}]}{\Gamma(\alpha)}_{1}\int^{1+\epsilon}
(\log\frac{1+\epsilon}{\tau})^{\alpha-1}f(\tau,x(\tau))\frac{d\tau}{\tau}\bigg]\Biggl|
\\&\quad\leq\frac{1}{\Gamma(\alpha)}_{1}\int^{t}(\log\frac{t}{\tau})^{\alpha-1}
\mid f(\tau,x(\tau))\mid\frac{d\tau}{\tau}~
\\&\quad\quad\quad\quad~+(\frac{\log t}{\log(1+\epsilon)})^{\gamma-1}
~\frac{1}{\Gamma(\alpha)}_{1}\int^{1+\epsilon}(\log\frac{1+\epsilon}{\tau})^{\alpha-1}
\mid f(\tau,x(\tau))\mid\frac{d\tau}{\tau}
\\&\quad\quad\quad+[\frac{1}{\log(1+\epsilon)}-\frac{1}{\log t}]\frac{(\log
t)^{\gamma-1}}{\sum_{i=0}^1\eta_{i}}\bigg[\frac{1}{\Gamma(\alpha-1)}_{1}\int^{e}(\log\frac{e}{\tau})^{\alpha-2}\mid
f(\tau,x(\tau))\mid\frac{d\tau}{\tau}
\\&\quad\quad\quad\quad\quad\quad+\frac{\nu}{\Gamma(\alpha-1)}_{1}\int^{\zeta}
(\log\frac{\zeta}{\tau})^{\alpha-2}\mid
f(\tau,x(\tau))\mid\frac{d\tau}{\tau}
\\&\quad\quad\quad+\frac{(1-\gamma)(\log(1+\epsilon))^{1-\gamma}[1-\nu(\log
\zeta)^{\gamma-2}]}{\Gamma(\alpha)}_{1}\int^{1+\epsilon}(\log\frac{1+\epsilon}{\tau})^{\alpha-1}\mid
f(\tau,x(\tau))\mid\frac{d\tau}{\tau}\bigg]
\\&\quad\leq C_{1}\bigg\{\frac{[1+(\log(1+\epsilon))^{1-\gamma+\alpha}]}{\Gamma(\alpha+1)}\\&\quad\quad\quad\quad\quad+
\frac{[1-(\log(1+\epsilon))]}{(\log(1+\epsilon))\Gamma(\alpha)\sum_{i=0}^1\eta_{i}}\bigg[1+\nu(\log
\zeta)^{\alpha-1}\\&\quad\quad\quad\quad+\frac{(1-\gamma)(\log(1+\epsilon))^{1-\gamma+\alpha}}{\alpha}[1-\nu(\log
\zeta)^{\gamma-2}]\bigg]\bigg\}=C_{2}\quad\quad\quad\quad\quad\quad\quad\quad\quad\quad\quad\quad\quad(3.10)
\end{align*}
this implies that $\parallel(\rho x)(t)\parallel\leq C_{2}.$
Moreover,
\begin{align*}
&\mid~ _{H}D^{1,1}(\rho x)(t)\mid\\&\quad=\Biggl|
-\frac{1}{\Gamma(\alpha-1)}_{1}\int^{t}(\log\frac{t}{\tau})^{\alpha-2}f(\tau,x(\tau))\frac{d\tau}{\tau}~
\\&\quad\quad\quad-\frac{(1-\gamma)}{\log(1+\epsilon)}(\frac{\log
t}{\log(1+\epsilon)})^{\gamma-2}~\frac{1}{\Gamma(\alpha)}_{1}\int^{1+\epsilon}
(\log\frac{1+\epsilon}{\tau})^{\alpha-1}f(\tau,x(\tau))\frac{d\tau}{\tau}
\\&\quad\quad\quad+[\frac{(\gamma-1)}{\log(1+\epsilon)}-\frac{(\gamma-2)}{\log t}]\frac{(\log
t)^{\gamma-2}}{\sum_{i=0}^1\eta_{i}}\bigg[\frac{1}{\Gamma(\alpha-1)}_{1}\int^{e}(\log\frac{e}{\tau})^{\alpha-2}f(\tau,x(\tau))\frac{d\tau}{\tau}
\\&\quad\quad\quad\quad\quad\quad-\frac{\nu}{\Gamma(\alpha-1)}_{1}\int^{\zeta}(\log\frac{\zeta}{\tau})^{\alpha-2}f(\tau,x(\tau))\frac{d\tau}{\tau}
\\&\quad\quad\quad+\frac{(1-\gamma)(\log(1+\epsilon))^{1-\gamma}[1-\nu(\log
\zeta)^{\gamma-2}]}{\Gamma(\alpha)}_{1}\int^{1+\epsilon}(\log\frac{1+\epsilon}{\tau})^{\alpha-1}f(\tau,x(\tau))\frac{d\tau}{\tau}\bigg]\Biggl|
\\&\quad\leq\frac{1}{\Gamma(\alpha-1)}_{1}\int^{t}(\log\frac{t}{\tau})^{\alpha-2}\mid f(\tau,x(\tau))\mid\frac{d\tau}{\tau}~
\\&\quad\quad\quad+\frac{(1-\gamma)(\log(1+\epsilon))^{1-\gamma}}{\Gamma(\alpha)}_{1}\int^{1+\epsilon}
(\log\frac{1+\epsilon}{\tau})^{\alpha-1}\mid
f(\tau,x(\tau))\mid\frac{d\tau}{\tau}
\\&\quad\quad\quad+\frac{(\gamma-1)[1-\log(1+\epsilon)]+\log(1+\epsilon)}{\log(1+\epsilon)
\sum_{i=0}^1\eta_{i}}\bigg[\frac{1}{\Gamma(\alpha-1)}_{1}\int^{e}(\log\frac{e}{\tau})^{\alpha-2}\mid
f(\tau,x(\tau))\mid\frac{d\tau}{\tau}
\\&\quad\quad\quad\quad\quad\quad+\frac{\nu}{\Gamma(\alpha-1)}_{1}\int^{\zeta}(\log\frac{\zeta}{\tau})^{\alpha-2}\mid
f(\tau,x(\tau))\mid\frac{d\tau}{\tau}
\\&\quad\quad\quad+\frac{(1-\gamma)(\log(1+\epsilon))^{1-\gamma}[1-\nu(\log
\zeta)^{\gamma-2}]}{\Gamma(\alpha)}_{1}\int^{1+\epsilon}(\log\frac{1+\epsilon}{\tau})^{\alpha-1}\mid
f(\tau,x(\tau))\mid\frac{d\tau}{\tau}\bigg]
\\&\quad\leq C_{1}\bigg\{\frac{1}{\Gamma(\alpha)}+\frac{(1-\gamma)(\log(1+\epsilon))^{1-\gamma+\alpha}}{\Gamma(\alpha+1)}
\\&\quad\quad\quad\quad\quad+\frac{(\gamma-1)[1-\log(1+\epsilon)]+\log(1+\epsilon)}{\log(1+\epsilon)\sum_{i=0}^1\eta_{i}}
\bigg[\frac{1+\nu(\log
\zeta)^{\alpha-1}}{\Gamma(\alpha)}\\&\quad\quad\quad\quad\quad\quad\quad\quad
+\frac{(1-\gamma)(\log(1+\epsilon))^{1-\gamma+\alpha}[1-\nu(\log
\zeta)^{\gamma-2}]}{\Gamma(\alpha+1)}\bigg]\bigg\}=C_{3}\quad\quad\quad\quad\quad\quad\quad\quad\quad(3.11)
\end{align*}
Thus,for each $t_{1},t_{2}\in J,$ we get
\begin{align*}
\mid(\rho x)(t_{1})-(\rho x)(t_{2})\mid  \leq
~~_{t_{1}}\int^{t_{2}}~_{H}D^{1,1}(\rho x)(\tau)\frac{d\tau}{\tau}
\leq
C_{3}(t_{2}-t_{1})\quad\quad\quad\quad\quad\quad\quad\quad\quad\quad\quad\quad(3.12)
\end{align*}
 which implies that $\rho$ is continuous over
$J$. Hence,the operator $~\rho:K\rightarrow K$ is completely
continuous,(by the Arzela-Ascoli theorem).
\\Finally,consider the set $U = \big\{v\in K\mid x
 =\lambda Tx, 0 < \lambda < 1\big\},$we show that the set $U$ is
bounded. Assume that $x\in U$,then $x=\lambda \rho x,\quad 0<
\lambda <1.$
\par Now for any $t\in J,$ we get
\begin{align*}
\mid x(t)\mid&=\lambda\mid(\rho x)(t)\mid
\\&\leq\frac{1}{\Gamma(\alpha)}_{1}\int^{t}(\log\frac{t}{\tau})^{\alpha-1}\mid
f(\tau,x(\tau))\mid\frac{d\tau}{\tau}
\\&\quad\quad\quad+(\frac{1}{\log(1+\epsilon)})^{\gamma-1}
\frac{1}{\Gamma(\alpha)}_{1}\int^{1+\epsilon}(\log\frac{1+\epsilon}{\tau})^{\alpha-1}\mid
f(\tau,x(\tau))\mid\frac{d\tau}{\tau}
\\&\quad\quad\quad\quad+\frac{[1-(\log(1+\epsilon))]}{(\log(1+\epsilon))\Gamma(\alpha-1)
\sum_{i=0}^1\eta_{i}}\bigg[_{1}\int^{e}(\log\frac{e}{\tau})^{\alpha-2}\mid
f(\tau,x(\tau))\mid\frac{d\tau}{\tau}
\\&\quad\quad\quad\quad\quad +{\nu}_{1}\int^{\zeta}(\log\frac{\zeta}{\tau})^{\alpha-2}\mid
f(\tau,x(\tau))\mid\frac{d\tau}{\tau}
\\&\quad\quad\quad+\frac{(1-\gamma)(\log(1+\epsilon))^{1-\gamma}[1-\nu(\log
\zeta)^{\gamma-2}]}{\alpha-1}~_{1}\int^{1+\epsilon}(\log\frac{1+\epsilon}{\tau})^{\alpha-1}\mid
f(\tau,x(\tau))\mid\frac{d\tau}{\tau}\bigg]
\\&\leq  C_{1}\bigg\{\frac{1}{\Gamma(\alpha+1)}
+\frac{(\log(1+\epsilon))^{1-\gamma+\alpha}}{\Gamma(\alpha+1)}
\\&\quad\quad\quad\quad+\frac{[1-(\log(1+\epsilon))]}{(\log(1+\epsilon))\Gamma(\alpha)\sum_{i=0}^1\eta_{i}}\bigg[1+\nu(\log
\zeta)^{\alpha-1}\\&\quad\quad\quad\quad\quad\quad\quad+\frac{(1-\gamma)(\log(1+\epsilon))^{1-\gamma+\alpha}}{\alpha}[1-\nu(\log
\zeta)^{\gamma-2}]\bigg]\bigg\}=M\quad\quad\quad\quad\quad\quad(3.13)
\end{align*}
\\Therefore, $\parallel x(t)\parallel \leq M$ for any $t\in
J.$ Hence, the set $U$ is bounded. So, from the above and by the
Theorem 2.11, the operator $\rho$ has at least one fixed point,that
implies to the problem$(1.1)$ has at least one solution.
\\ \textbf{\ Theorem 3.6}  Assume that there exist a small positive
number $\tilde{r}$ and $0<\mu<\frac{1}{\Phi}$ such that
\[\mid f(t,x)\mid\leq\mu \mid x\mid ~~for~ 0<\mid x\mid < \tilde{r}, ~~ where ~~\Phi ~~is ~~defined ~~by ~(3.7).\]
Then the problem (1.1) has at least one solution.
\\\textbf{\ Proof.}  Firstly, let $K$ be a Banach space defined by
(3.5), and define \\$B_{\tilde{r}}=\{x\in K:\parallel
x\parallel\leq\tilde{r}\}$ and put $x\in K$ such that $\parallel
x\parallel=\tilde{r},$ that is, $x\in\partial B_{\tilde{r}}.$
\par Now, with the same argument of proof in the previous theorem, we
can shown that $\rho$ is completely continuous and we have,
\begin{align*}
&\parallel(\rho x)(t)\parallel\\&\quad=sup_{t\in J}\bigg\{\Biggl|
-\frac{1}{\Gamma(\alpha)}_{1}\int^{t}(\log\frac{t}{\tau})^{\alpha-1}f(\tau,x(\tau))\frac{d\tau}{\tau}~
\\&\quad\quad\quad\quad\quad\quad\quad+(\frac{\log t}{\log(1+\epsilon)})^{\gamma-1}~\frac{1}{\Gamma(\alpha)}_{1}\int^{1+\epsilon}
(\log\frac{1+\epsilon}{\tau})^{\alpha-1}f(\tau,x(\tau))\frac{d\tau}{\tau}
\\&\quad\quad\quad+[\frac{1}{\log(1+\epsilon)}-\frac{1}{\log t}]\frac{(\log
t)^{\gamma-1}}{\sum_{i=0}^1\eta_{i}}\bigg[\frac{1}{\Gamma(\alpha-1)}_{1}\int^{e}
(\log\frac{e}{\tau})^{\alpha-2}f(\tau,x(\tau))\frac{d\tau}{\tau}
\\&\quad\quad\quad\quad\quad\quad-\frac{\nu}{\Gamma(\alpha-1)}_{1}\int^{\zeta}
(\log\frac{\zeta}{\tau})^{\alpha-2}f(\tau,x(\tau))\frac{d\tau}{\tau}
\\&\quad\quad\quad+\frac{(1-\gamma)(\log(1+\epsilon))^{1-\gamma}[1-\nu(\log
\zeta)^{\gamma-2}]}{\Gamma(\alpha)}_{1}\int^{1+\epsilon}(\log\frac{1+\epsilon}{\tau})^{\alpha-1}
f(\tau,x(\tau))\frac{d\tau}{\tau}\bigg]\Biggl|\bigg\}\\&\quad\leq\mu\Phi\parallel
x\parallel\quad\quad\quad\quad\quad\quad\quad\quad\quad\quad\quad\quad\quad\quad
\quad\quad\quad\quad\quad\quad\quad\quad\quad\quad\quad\quad\quad\quad\quad\quad\quad\quad(3.14)
\end{align*}
\\Hence, $\parallel(\rho x)(t)\parallel\leq\parallel x\parallel,$
with $x\in\partial B_{\tilde{r}}.$ Then by applied the theorem 2.12,
\\$\rho$ has at least one fixed point. Therefore, the problem(1.1) has
at least one solution on $J.$
\\ \textbf{\ Example.} \\ Consider the following boundary value
problem for Hilfer-Hadamard-type fractional differential equation:
$$\quad\quad\quad\quad\quad\quad\quad\quad_{H}D^{3/2,1/2}x(t)+f(t,x(t))=0,~\quad\quad~~t\in
J:=(1,e]~\quad\quad\quad\quad\quad\quad\quad\quad\quad\quad(3.15)$$
$$\quad x(1.2)=0,~\quad\quad~_{H}D^{1,1}x(e)=(1/2)~_{H}D^{1,1}x(3/2).~$$
Here,
$$\alpha=3/2,\quad\beta=1/2,\quad\gamma=7/4,\quad\nu=1/2,\quad\zeta=3/2,\quad\epsilon=0.2,\quad1+\epsilon=1.2
$$ and $$f(t,x(t))=\frac{1}{32}(\sqrt{t}+\log t)\big(\frac{\mid x\mid}{2+\mid x\mid}\big). $$
Clearly, $$\mid f(t,x)\mid\leq\frac{1}{32}(\sqrt{t}+1 )(\mid
x\mid+1)$$ and $$\mid f(t,x)-f(t,y)\mid\leq\frac{1}{32}(\sqrt{t}+1
)(\mid x-y\mid)\leq\frac{1}{16}\mid x-y\mid$$ Therefore, by
Theorem3.4, the boundary value problem (3.15) has a unique solution
on $(1,e]$ with $C=\frac{1}{16}=0.0625.$ We can show that
$\Phi=1.404,\quad C\Phi=0.0876<1.$

\end{document}